\documentstyle[11pt]{article}

\font\sit=cmti8

\font\smc=cmcsc10

\def\uf{{\cal U}_{\varphi}}

\def\f{\varphi}

\def\Mf{M^\varphi}

\def\ben{\begin{enumerate}}
\def\een{\end{enumerate}}

\def\zZ{{\rm Z}\!\!{\rm Z}}
\def\zC{{\rm C}\!\!\!\vrule height 6.4pt depth -0.4pt width
1pt \ }
\def\zN{{{\rm I}\!{\rm N}}}

\newtheorem{lemm}{Lemma}[section]
\newtheorem{theo}[lemm]{Theorem}
\newtheorem{prop}[lemm]{Proposition}
\newtheorem{coro}[lemm]{Corollary} 
\newtheorem{rema}[lemm]{Remark} 
 
\newtheorem{exam}[lemm]{Example}
\newtheorem{exams}[lemm]{Examples}

\newenvironment{proof}[1]{
  \trivlist \item[\hskip \labelsep{\bf #1}]}{\hfill\mbox{$\Box$}
  \endtrivlist}

\begin{document}

\centerline{\bf \Large HOMOTOPY OF STATE ORBITS}
\vskip.6cm

\centerline{\large \smc Esteban Andruchow and Alejandro Varela}
\vskip.7cm  
\centerline{\large Instituto de Ciencias -- Univ. Nac. de Gral. Sarmiento}
\centerline{\large Argentina}
\vskip1cm  


\vskip0.5cm

\begin{abstract} 
\noindent 
Let $M$ be a von Neumann algebra, $\f$ a 
faithful normal state and denote by $M^\f$ the fixed point 
algebra of the modular group of $\f$. Let $U_M$ and $U_{M^\f}$ be the unitary 
groups of $M$ and $M^\f$. In this paper we study the quotient
$\uf =U_M/U_{M^\f}$ endowed with two natural topologies: the one induced by the usual 
norm of $M$ (called here usual topology of $\uf$), and the one induced by the pre-Hilbert
C$^*$-module norm given by the $\f$-invariant conditional expectation 
$E_{\f}:M \to M^{\f}$ (called the modular topology). It is shown that $\uf$ 
is simply connected with the usual topology. Both topologies are compared, and 
it is shown that they coincide if and only if the Jones index of $E_{\f}$ is 
finite. The set $\uf$ can be regarded as a model for the unitary orbit
$\{\f \circ Ad(u^*): u\in U_M\}$ of $\f$, and either with the usual or the 
modular it can be embedded continuously in the conjugate space $M^*$ (although
not as a topological submanifold).

\end{abstract}

\section{Introduction}
Let $M$ be a von Neumann algebra and $\f$ a faithful normal state of $M$. Denote by 
$U_M$ the unitary group of $M$, and by $M^\f$ the centralizer of $\f$, that is 
$M^\f=\{x\in M : \f(xy)=\f(yx) \hbox{ for all } y\in M\}$. Let $\uf$ be  the unitary
orbit of $\f$, i.e.
$$
{\cal U}_\f=\{ \f \circ Ad(u) : u\in U_M\}
$$
where $Ad(u)(x)=uxu^*$. 
The isotropy subgroup at $\f$ (=the set of unitaries
that leave $\f$ fixed) is the unitary group $U_{M^\f}$ of $M^\f$. 
In previous papers we introduced a homogeneous and reductive structure 
for $\uf$, by means of the natural identification
$$
\uf \simeq U_M / U_{M^\f} .
$$
We are not regarding $\uf$ with the norm
topology of $M^*$ (in \cite{av2} it was shown that in general $\uf$ is not a submanifold 
of $M^*$), but with the quotient topology induced by the usual norm of $M$. With this
topology $\uf$ is a real analytic manifold. 
In particular, the map
$$
\pi_\f : U_M \to \uf \  ,  \ \pi_\f (u)=\f\circ Ad(u^*)
$$
is a (principal) fibre bundle,
with fibre $U_{M^\f}$.
In  section 2 we use this fibration to prove that these orbits $\uf$ are always simply connected.

There is another natural topology in the set $U_M / U_{M^\f}$. Namely, let 
$E=E_\f$ be the unique $\f$-invariant conditional expectation $E:M \to M^\f$. 
This gives rise to a natural pre-Hilbert C$^*$-module structure for $M$, with $M^\f$ 
valued inner product given by $<x,y>=E(x^*y)$ and norm $\|x\|_E=\|E(x^*x)\|^{1/2}$. It 
is well known that $M$ is $\| \ \|_E$-complete if and only if the Jones index of $E$ is 
finite (\cite{bdh},\cite{pp}). The condition that $E$ be of finite index is a rather strict
 requirement for  $\f$. This implies that in general the usual norm and $\| \ \|_E$ 
  define different topologies for $M$ and for $U_M / U_{M^\f}$. We  call them, respectively, 
  the {\bf usual} and the {\bf modular} topology in $\uf$. 

One has the following inequality, for $u,w \in U_M$:
$$
\| \f \circ Ad(u^*) - \f \circ Ad(w^*)\| \le 2\|u - w\|_E \le 2\|u - w\|
$$
where the first norm is the usual norm in $M^*$. If one replaces $u, w$ by 
$uv, wv'$, for $v,v' \in U_{M^\f}$ then $\f \circ Ad((uv)^*)=\f \circ Ad(u^*)$ 
and $\f \circ Ad((wv')^*)=\f \circ Ad(w^*)$. This implies that $\uf$, both in 
the usual and the modular topology, can be embedded in the conjugate space $M^*$. 
These matters are discussed in section 3. We present different models for $\uf$, 
inside the grassmannians of the basic extension of $E$ (see below), and inside the 
interior tensor product $M \otimes_{M^\f} M$ of $M$ regarded as a pre-Hilbert module $M$.

Finally, let us recall the basic extension of $E:M \to M^\f$. Denote 
by $H_\f$ the completion of the pre-Hilbert space $M$ with the inner 
product given by $\f$. Then $E$ is bounded for this inner product, and 
therefore extends to a selfadjoint projection $e=e_\f$, called the Jones 
projection, whose range is the closure of $M^\f$ in $H_\f$. Let $M_1\subset B(H_\f)$ 
be the von Neumann algebra generated by $M$ and $e$. We refer the reader to 
\cite{jones}, \cite{pp} or \cite{bdh} for the details of this construction. 
Some of the properties of $e$ are:
\begin{enumerate}
\item[-] $eae=E(a)e, \ a\in M$
\item[-] $M\cap \{e\}'=M^\f$
\item[-] The map $x \mapsto xe$ is a *-isomorphism between $M^\f$ and $M^\f e$.
\end{enumerate}
\noindent

\bigskip

\section{The fundamental group of $\uf$ is trivial}
Throughout this section $\uf$ is endowed with the usual topology (i.e. 
the quotient topology induced by the usual norm of $M$).
Recall that a fibre bundle gives rise to an exact sequence of homotopy 
groups. In our case, the bundle $\pi_\f$ yields the exact sequence
$$
\cdots \pi_2 (\uf) \to \pi_1(U_{M^\f}) \stackrel{i_*}{\to}
\pi_1 (U_M) \to \pi_1 ({\cal U_\f}) \to \pi_0 (U_{M^\f}) = 0 ,
$$
where $1$ is taken as base point for the homotopy groups of the unitary 
groups and $\f$ is the base point for $\uf$. Here $i_*$ denotes the homomorphism 
induced by the inclusion $i : U_{M^\f} \hookrightarrow U_M$.
We can use then results by Handelman \cite{han}, Schr\"{o}der \cite{sch}, 
Breuer \cite{bre}, as well as the classical result by Kuiper \cite{kui}, computing the homotopy groups of
the unitary group of a von Neumann algebra, in order to 
obtain information about $\uf$. 

Note that the center ${\cal Z}(M^\f)$ of $M^\f$ includes 
the center ${\cal Z}(M)$ of $M$. Suppose that $p\in {\cal Z}(M)$
is a projection. Then $\f _p=\f|_{Mp}$ is a faithful and
normal state in $Mp$ whose centralizer is the algebra $M^\f  p$,
and the canonical conditional expectation $E_{\f _p}$ is the restriction
of $E$ to $Mp$. In other words, each  projection in the center of $M$
factorizes the unitary groups of $M$ and $M^\f$ and the orbit $\uf$:
$$
U_M \simeq U_{Mp} \times U_{M(1-p)} \ , \ U_{M^\f} \simeq U_{(Mp)^{\f _p}}
\times U_{(M(1-p))^{\f _{1-p}}} \ \hbox{ and } \uf \simeq {\cal U}_{\f _p}
\times {\cal U}_{\f _{1-p}}.
$$

Therefore, if one considers the
central type decomposition projections of $M$, the study 
of the homotopy group of $\uf$ reduces to the case when $M$
is of a definite type. We shall proceed to show that $\uf$ 
is simply connected through a series of lemmas, covering the 
possible types of $M$ and $M^\f$.

By well known results \cite{kui} \cite{bre}, the properly 
infinite part of $M$ gives state orbits with vanishing $\pi_1$ group. Indeed,

\begin{lemm}
If $M$ is a properly infinite von Neumann algebra, then $\uf$ is simply connected.
\end{lemm}
\begin{proof}{Proof.}
Since $\uf$ is connected, it remains to prove that $\pi_1(\uf)$ is trivial. 
This follows by appealing to the homotopy exact sequence, and the fact 
(\cite{kui}, \cite{bre}) that $U_M$ has trivial $\pi_1$ group.
\end{proof}

One is therefore constrained to the case when $M$ is finite. 
Let us recall the following results 
(see \cite{acs}), which are based on the results of \cite{han} 
(see also \cite{sch}). If $M$ is of type II$_1$, 
Handelman proved that $\pi_1(U_M)$ is isomorphic to (the additive group) 
${\cal Z}(M)_{sa}$ of selfadjoint elements of the center of $M$.

\begin{lemm}
Let $M$ be a type II$_1$ von Neumann algebra and $Tr$ its center valued trace. 
Suppose that $N \subset M$ is a von Neumann subalgebra with the same unit, and 
denote by $i$ the inclusion map $i: U_N \hookrightarrow U_M$. Then the image
of  the homomorphism $i_* :\pi_1 (U_N) \to \pi_1 (U_M)\simeq {\cal Z}(M)_{sa}$ 
is equal to the additive group generated by
the set $\{Tr(p) : p \hbox{ projection in } N\}$.
\end{lemm}

There is an analogous result for the type I case. If $M$ is of type 
I$_n$, then $\pi_1(U_M)\simeq C(\Omega, \zZ)$, where $\Omega$ denotes the Stone space of the center of $M$. 

\begin{lemm}
Let $M$ be a von Neumann algebra of type I$_n$, $N\subset M$, and $Tr$  
the center valued trace of $M$, then the image of $i_*$ identifies with the group generated by the functions 
$\{Tr(p)\hat{ }: p \hbox{ projection in } N\}$.
\end{lemm}

In order to compute the $\pi_1$ group of $\uf$ we shall apply these results to the 
case $N=M^\f$. From these lemmas it is clear that one needs to compute 
$\{Tr(p) : p \hbox{ projection in } M^\f\}$. 

The next step is to try further reductions using the type decomposition central projections of $M^\f$.

\begin{lemm}\label{2-2}
Suppose that $M$ and $M^\f$ are of type II$_1$. Then $i_*$ is surjective. 
As a consequence,  $\uf$ is simply connected.    
\end{lemm}
\begin{proof}{Proof.}
We claim that in this case the homomorphism $i_*$ identifies with 
$Tr|_{{\cal Z}(M^\f)_{sa}} :{\cal Z}(M^\f)_{sa} \to {\cal Z}(M)_{sa}$. 
Then it is clear that $i_*$ is surjective. In order to prove our claim, 
we shall see first that if $Tr$ and $Tr^\f$ denote respectively the center 
valued traces of $M$ and $M^\f$, then $Tr\circ Tr^\f=Tr|_{\Mf}$. Indeed, let $x\in \Mf$, then
$Tr^\f (x)$ is the (norm) limit of a sequence of elements $v_nxv_n^*$ in 
$co\{ vxv^*: v\in U_{\Mf}\}\cap {\cal Z}(\Mf)$ (where $co\{ vxv^*: v\in U_{\Mf}\}$
 denotes the convex hull of $\{ vxv^*: v\in U_{\Mf}\}$). Since $Tr(v_nxv_n^*)=Tr(x)$, 
 it follows that $Tr(Tr^\f(x))=Tr(x)$.
Now let $z$ be an element in ${\cal Z}(\Mf)$ with $0\le z \le 1$, then there exists 
 a projection $p \in \Mf$ such that $Tr^\f (p)=z$. Under the identification 
 $\pi_1(U_{\Mf})\simeq {\cal Z}(\Mf)_{sa}$, 
the class of the loop $\gamma (t)=e^{itp}$ in $U_{\Mf}$ corresponds to the element $z$,
$i_*$ sends this element to the class of $\gamma$ in $U_M$, that is, 
to $Tr(p)$ (see \cite{han}). Therefore
$i_*(z)=Tr(p)=Tr(Tr^\f (p))=Tr(z)$. The fact that $\uf$ is simply connected 
follows from the exact sequence, where $\pi_1(\uf)=\pi_1(U_M)/Im(i_*)=0$.
\end{proof}
Let us consider the following examples, which show that for $M$ of type II$_1$, 
different types of $\Mf$ can occur.
\begin{exams}

\begin{enumerate}

\item[1]
Suppose that $M$ is of type II$_1$, let $p$ be a projection and $\tau$ a faithful tracial state. 
Consider $h=\frac{1}{2s}p +\frac{1}{2(1-s)}(1-p)$, with $s=\tau(p)$. Put $\f(x)=\tau(hx)$. 
Clearly $\f$ is a faithful and normal state, with $E (x)=pxp +(1-p)x(1-p)$ and 
$\Mf=\{p\}'\cap M=pMp \oplus (1-p)M(1-p)$, which is also of type II$_1$. 
A similar example can be done with a family $\{p_n \in M\}_{n\in \zN}$ of mutually orthogonal projections.

\item[2]
Let again $M$ be of type II$_1$, and let ${\cal A} \subset M$ be a maximal abelian 
subalgebra. Choose $\tau$ a faithful normal tracial state and $h$ a positive operator 
without kernel (i.e. $ha=0$ implies $a=0$) that generates ${\cal A}$, normalized so 
that $\tau (h)=1$. Note that ${\cal A}$ is purely non atomic. Consider the state 
$\f(x)=\tau(hx)$. Then $\Mf=\{h\}'\cap M ={\cal A}$.
Another example can be obtained by tensoring $M$ with $M_n(\zC)$ and $\f$ with the 
usual trace $t_n$ of $M_n(\zC)$. In this case $M_1=M\otimes M_n(\zC)$ is of type II$_1$ 
and $M_1^{\f \otimes t} =\Mf \otimes M_n(\zC)^{t_n}=M_n({\cal A})$ which is of type I$_n$. 
\end{enumerate}
\end{exams}
Note that if $M$ is a finite von Neumann algebra with a faithful normal state $\f$ 
such that $\Mf$ is abelian, then $\Mf$ must be maximal abelian in $M$. Indeed, 
if $a\in (\Mf) ' \cap M$, and $\f=\tau (h \ )$ for a tracial state $\tau$, and $x\in M$, 
then $\f (ax)=\tau (hax)=\tau (ahx)=\tau (hxa)=\f(xa)$,
i.e. $a \in \Mf$. In particular, if $M$ is of type II$_1$, 
then the (necessarily non trivial) center of $\Mf$ must be non atomic, and therefore
the situation in the above example, part 2, is essentially the only possible one.

\begin{lemm}
If $M$ is of type II$_1$ and $\Mf$ is abelian, then $\uf$ is simply connected.
\end{lemm}

\begin{proof} {Proof.}
As before, one needs to show that if $Tr$ is the center valued trace of $M$, then 
$\{ Tr(q): q \hbox{ projection in } \Mf \}={\cal Z}(M)_{sa}$. First, pick $c \in {\cal Z}(M)_{sa}$ 
of the form $c=\sum_{i=1}^{n} \alpha_i p_i$
with $p_i$ mutually orthogonal and $0\le \alpha_i \le 1$. Given $\epsilon >0$ sufficiently small, 
we claim that there exists projections $q_i \in \Mf$ such that 
$0\le (\alpha_i - \epsilon /n)p_i \le Tr(q_i) \le \alpha_i p_i$.
Indeed, otherwise there would be a central projection $p$ and an 
interval $(0, \lambda)$ such that between $0$ and $\lambda p$ there 
are no values $Tr(q)$ with $q$ projection in $\Mf$. In that case, put 
$$
\lambda_0 = \sup \{\lambda>0 : \hbox{ there are no projections } q\in \Mf \hbox{ with } Tr(q) \le \lambda p \}.
$$
Clearly, $0<\lambda_0\le 1$. Then one can find sequences $\lambda_n > \lambda_0$ and $q_n$, 
where  $\lambda_n$ decreases to $\lambda_0$ and $q_n$ are projections in $\Mf$ with $q_n \ge q_{n+1}$, 
such that $Tr(q_n)\le \lambda_n p$. Then $q_0=\wedge_{n} q_n$ is a projection in $\Mf$ 
satisfying $Tr(q_0)=\lambda_0 p$. Suppose now that there exists a projection $r$  in $\Mf$ with $0< r
< q_0$. Let $e_n = \chi_{[0,\lambda_0- 1/n)}(Tr(r))$ be the spectral projection of $Tr(r)$ 
associated to the interval $[0,\lambda_0 -1/n)$, lying in ${\cal Z}(M)$. Then clearly $r_n=e_n r$ 
is a projection in $\Mf$ satisfying that $Tr(r_n) \le \lambda_0 - 1/n$. Since $r_n$ increases 
to $r$, there exists $n$ such that $r_n $ is non zero, this implies a contradiction with the 
fact that $\lambda$ is the supremum of the above set. Therefore no
such $r$ should exist, which in turn would imply that $q_0$ is a minimal projection in $\Mf$. 
Since $\Mf$ is maximal abelian in $M$ of type II$_1$, it has no minimal projections, and we 
arrive to a contradiction. Returning to our original central element $c$, it follows that 
we can find projections $q_i$ in $\Mf$ with 
$0\le (\alpha_i - \epsilon /n)p_i \le Tr(q_i) \le \alpha_i p_i$. Since $q_i \le p_i$ and 
these projections are mutually orthogonal, it follows that $q=\sum_{i=1}^n q_i$ is a projection in $\Mf$.
Moreover, we have that $c-\epsilon \le Tr(q) \le c$. Then we can construct an increasing 
sequence $q_n$ of projections in $\Mf$ such that $Tr(q_n)$ converges to $c$. Pick $q'=\vee q_n$, clearly $Tr(q')=c$.
Now, if $c$ is any element in ${\cal Z}(M)$ with $0\le c \le 1$, let $c_n$ be an increasing sequence 
of positive elements in ${\cal Z}(M)$ with finite spectrum, converging to $c$ in norm. We can find 
projections $q'_n$ in $\Mf$ with $Tr(q'_n)=c_n$. Then $q'_n$ is an increasing sequence of projections, 
put $q'_0=\vee q'_n$. We obtain that $Tr(q'_0)=c$, and the proof is complete.

\end{proof}

\begin{lemm}\label{2-1}
Let $M$ be a von Neumann algebra of type II$_1$ and $\f$ a faithful and normal state 
such that $\Mf$ is of type $I$. Then $\uf$ is simply connected. 
\end{lemm}
 \begin{proof}{Proof.}
Let $p_n$ be the projections of the center of $\Mf$ decomposing it in its type I$_n$ parts, 
$n<\infty$. Pick $c\in {\cal Z}(M)$, and put $c_n=cp_n$. Suppose that for each $n$ we can 
find $q_n$ in $p_n \Mf\subset \Mf$ with $Tr(q_n)=c_n$. Then $q=\sum_n q_n$ is a projection 
in $\Mf$ such that $Tr(q)=c$. Therefore
it remains to prove our statement in the case $\Mf$ of type I$_n$. Indeed,
note that $p_n \Mf$ is the centralizer of the (faithful and normal)
state $\f_n$ of $p_n M p_n$,
which is simply the restriction of $\f$ to $p_n Mp_n$.
 
Let now $e$ be a minimal abelian  projection in $\Mf$. Again, pick $0\le c \le 1$ in 
${\cal Z}(M)$. Now $e M e$ is of type II$_1$, and the state $\f_e$ of $e M e$ given 
by the restriction of $\f$ to this algebra has centralizer equal to $e \Mf e$. By the lemma above, there
exists a projection $q \in e \Mf e \subset \Mf$ such that 
$$
Tr_e(q)=ec,
$$ 
where $Tr_e$ is the center valued trace of $e M e$, i.e.
$Tr_e(e x e)=e Tr(e x e)$. Since $\Mf$ is of type I$_n$, it follows that
$Tr (e)=1/n$. Taking trace in the above equality yields
$(1/n) Tr(q)=(1/n) c$ and the statement follows.
\end{proof}
Finally, the case when $M$ is of type $I$ is dealt in a similar way.

\begin{lemm}\label{1-1}
If $M$ is a finite type I von Neumann algebra, and $\f$ is a faithful and normal state, then $\uf$ is  simply connected.
\end{lemm}
\begin{proof}{Proof.}
As remarked before, one can restrict to the case when $M$ is of type I$_n$, for a fixed $n<\infty$. 
In this case, $\pi_1(U_M)$ equals $C(\Omega,\zZ)$, where the isomorphism is implemented by the map 
sending the class of the curve $\alpha (t)= e^{2\pi i t p}$ to the continuous map $Tr(p)\hat{ }$, 
for $p$ a projection in $M$. In other words, $\pi_1(U_M)$ identifies with elements $c \in {\cal Z}(M)_{sa}$ 
which are of the form $c=\sum_{i=1}^k m_i p_i$, with $p_i$ mutually orthogonal in ${\cal Z}(M)$ and $m_i$ 
are integers. The proof follows, recalling that ${\cal  Z}(M) \subset {\cal Z}(M^\f)$, and therefore 
$Tr(p_i)=p_i$, i.e. $c$ lies in the image of $i_*$.
\end{proof}
We may state now our theorem:
\begin{theo}
Let $M$ be a von Neumann algebra, and $\f$ a faithful and normal state. The  unitary orbit of $\f$, 
$\uf =\{ \f \circ Ad(u): u \in U_M \}$ regarded with the (usual) quotient topology $U_M/U_{\Mf}$ is simply connected.
\end{theo}
\begin{proof}{Proof.}
As noted at the beginning of the section, it suffices to prove the statement in the case when $M$ is 
of a definite (finite) type. Type I case was dealt in \ref{1-1}. Suppose that $M$ is of type II$_1$.
Then $M^\f$ is finite, and there exist two projections $p_I, p_{II}$ in 
${\cal Z}(M^\f)$ such that $p_I+p_{II}=1$, $p_IM^\f$ is of type I and 
$p_{II}M^\f$ is of type II$_1$. In \cite{acs2} it was shown that if $p$ 
is a projection in a von Neumann algebra $M$, then the unitary orbit 
$\{upu^*: u\in U_M\}$ is simply connected. This unitary orbit is the base space 
of a fibration of $U_M$ with fibre $U_N$ where $N=\{p\}'\cap M=\{pxp +(1-p)x(1-p): x\in M\}$. 
In other words, the quotient $U_M/U_N$ is simply connected. In 
our case we have that 
$$
U_M/\left(U_{p_IMp_I }\times U_{p_{II}Mp_{II}}\right)
$$
is simply connected. The inclusion $U_{M^\f} \subset U_M$ can be factorized
$$
U_{M^\f}=U_{p_IM^\f}\times U_{p_{II}M^\f} \subset U_{p_IMp_I} \times U_{p_{II}Mp_{II}} \subset U_M .
$$
The inclusion $U_{p_IM^\f} \subset  U_{p_IMp_I}$ induces an epimorphism of the $\pi_1$ groups, 
by lemma \ref{2-1}. The same happens with the inclusion $U_{p_{II}M^\f} \subset U_{p_{II}Mp_{II}}$, 
by lemma \ref{2-2}. The last inclusion $U_{p_IMp_I} \times U_{p_{II}Mp_{II}} \subset U_M$ also induces 
an epimorphism of the $\pi_1$ groups by the remark above. Therefore $\uf$ is simply connected also in this case.
\end{proof}

\section{Topologies in $\uf$}
In the previous section we considered in $\uf$ the quotient topology $U_M/U_{\Mf}$ with 
$U_{\Mf} \subset U_M$ endowed with the norm topology of $M$. In this section we shall 
consider in $M$ also the norm $\| \ \|_{E}$ given by $\|x\|_{E}=\|E(x^*x)\|^{1/2}$. 
That is, the norm of $M$ regarded as a pre-C$^*$-module over $\Mf$, with the $\Mf$ valued 
inner product $<x,y>=E(x^*y)$. It is known \cite{bdh} that $M$ is complete with this norm 
if and only if the index of $E$ is finite. We shall denote
the corresponding topologies induced in $\uf\simeq U_M/U_{\Mf}$ as the usual topology 
and the modular topology. Let us recall some facts

\begin{rema}
If the index of $E$ is finite, then both topologies coincide, because
both norms are equivalent in $M$ in this case. 

In general  (\cite{av2}) $\uf$ with the modular topology is naturally homeomorphic to the orbit
$$
U_M(e)=\{ue u^* : u \in U_M\}\subset M_1
$$
via the map $\f\circ Ad(u^*) \mapsto ueu^*$.
Here  $U_M(e)$ is considered with the norm topology of $M_1$. This orbit is a subset of 
the grassmannians (=projections) of $M_1$. It is a submanifold of the grassmannians if 
and only if the index of $E$ is finite. It follows that $\uf$ is in general a metric space, 
and a complete metric space in the finite index case.
\end{rema}

\begin{rema}
The condition that the centralizer expectation $E$ of a state $\f$ be of finite index is rather strong. 
It implies that $M$ must be finite. Moreover, if $M$ is a factor, it happens if and only if the Radon-Nikodym derivative of
$\f$ with respect to the trace of $M$ is a (bounded) operator with finite spectrum (see \cite{av2}). 
If this condition holds, then $\uf$ is simply connected with the modular topology as well.
\end{rema}

The following results establish that the finite index case is the only situation in which both topologies in $\uf$ coincide.
\begin{prop}
Let $F:M \to N\subset M$ be a faithful conditional expectation of infinite index. 
Then the norm of M and the norm $\| \ \|_F$ induced by $F$ define topologies in $U_M/U_N$ which are not equivalent.
\end{prop}
\begin{proof}{Proof.}
Since the index of $F$ is infinite \cite{bdh}, \cite{fk}, there exist elements $a_n \in M$  
with $0\le a_n \le 1$, $\|a_n\|=1$ and $F(a_n)\to 0$ as $n$ tends to infinity.  
It is straightforward to verify that the distance $d(a_n,N)=\inf \{\|a_n -b\|: b\in N\}$ 
does not tend to zero with $n$.
Let $u_n\in U_M$ be unitaries such that $1-a_n=\frac{u_n+u_n^*}{2}$. 
Then
$$
\|u_n -1\|_F^2=\|2-F(u_n +u_n^*)\|=2\|F(a_n)\|\to 0.
$$
Therefore the sequence of the classes of the elements $u_n$ tends to the class of $1$ 
in the modular topology. We claim that $[u_n]$ does not tend to $[1]$ in the usual 
topology (induced by the norm of $M$). Suppose not. Then there exist unitaries $v_n \in U_N$ 
such that $u_nv_n \to 1$. Then 
$$
\|u_n -v_n^*\|^2=\|(u_n-v_n^*)(u_n^*-v_n)\|=\|2-u_nv_n-v_n^*u_n^*\|\to 0.
$$
This implies that $d(u_n,N)\to 0$, and therefore $d(a_n,N)\to 0$, an absurd.
\end{proof}
\begin{coro}
The usual and the modular topology coincide in $\uf$ if and only if the index of $E$ is finite.
\end{coro}
In \cite{as2} it was shown that when the index of a conditional expectation $F:M \to N$ 
is finite then the mapping $U_M \ni u \mapsto ufu^* \in U_M(f)=\{ufu^*: u\in U_M\}$ 
is a (principal) fibre bundle (where $f$ denotes the Jones projection of $F$). Using 
the result above it can be shown that also the converse is true:
\begin{coro}
Let $F:M \to N$ be a conditional expectation and $f$ the Jones projection of $F$. 
Then the mapping
$$
U_M \ni u \mapsto ufu^* \in U_M(f)=\{ufu^*: u\in U_M\}
$$
has continuous local cross sections if and only if the index of $F$ is finite.
\end{coro}
\begin{proof}{Proof.}
It only remains to prove that if the above mapping has local cross sections, then 
the index of $F$ is finite. The existence of local cross sections implies that the 
bijective and continuous map induced in the quotient,
$$
U_M/U_N \to U_M(f)
$$
is open, and therefore a homeomorphism. On the other hand it holds in general 
(\cite{av2}) that this same bijection is a homeomorphism between $U_M(f)$ 
and the modular topology in $U_M/U_N$. It follows by the proposition
above, that the index of $F$ is finite.
\end{proof}

Next we show that $\uf$ with the modular topology, 
can be presented as a subset of the interior tensor product 
$M \otimes_{\Mf} M$ of the pre-C$^*$-module $M$ over $\Mf$ 
with itself (see  \cite{la} for the particulars of this 
construction). The inner product and the norm of 
$M \otimes_{\Mf} M$ are given by:
if $x_i, y_i \in M$, $i=1,2$ then
$$
<x_1\otimes y_1, x_2\otimes y_2>= E (y_1^*E(x_1^*x_2)y_2)
$$
and
$$
\|x_1\otimes y_1\|=\|E(y_1^*E(x_1^*x_1)y_1)\|^{1/2}.
$$
Consider the set 
$$
{\cal D}_\f=\{ u\otimes u^* : u \in U_M\}\subset M \otimes_{\Mf} M
$$
Clearly the map $U_M \to {\cal D}_\f$, $u\mapsto u\otimes u^*$ induces a well 
defined bijection
$$
\delta: U_M/U_{\Mf} \to {\cal D}_\f 
\qquad
\delta([u])=u\otimes u^*.
$$
\begin{prop}
The map $\delta$ is continuous both in the usual and modular topology of $\uf$. It is a homeomorphism in the modular topology.
\end{prop}
\begin{proof}{Proof.}
Since the space considered is homogeneous and the action of the unitary group is continuous, 
it suffices to consider continuity at the class $[1]$ of $1$.
First, note that  $u_\alpha \otimes u_\alpha^* \to 1\otimes 1$ in the norm topology of 
$M \otimes_{\Mf} M$ if and only if $E (u_\alpha)E(u_\alpha^*) \to 1$ in $\Mf$.
Then it is clear that the map $U_M \ni u \mapsto u\otimes u^* \in  M \otimes_{\Mf} M$ 
is continuous in the norm topology of $M$.
Therefore the map induced in the quotient $U_M/U_{\Mf}$, i.e. $\delta$,
is continuous in what we are calling the usual topology of the quotient.

In the modular topology, as noted above, ${\cal U}_\f$
is homeomorphic to the orbit $U_M(e)=\{ue u^*: u \in U_M\}\subset M_1$
in the norm topology. Therefore $[u_\alpha]\to [1]$ if and only if
$u_\alpha e u_\alpha^* \to e$.
This implies that $e u_\alpha e u_\alpha^* e=
E(u_\alpha)E_\alpha (u_\alpha^*)e \to e$.
Since the mapping $\Mf \to \Mf e$, $x\mapsto x e$ is a *-isomorphism,
it follows that $E (u_\alpha)E(u_\alpha^*)\to 1$, that is
 $u_\alpha \otimes u_\alpha^* \to 1\otimes 1$ in $M \otimes_{\Mf} M$.

In order to see that $\delta$ is a homeomorphism with the modular 
topology, suppose that $u_\alpha \in U_M$ such that $E (u_\alpha)E(u_\alpha^*)\to 1$. 
Therefore there exists $\alpha_0$ such that for $\alpha \ge \alpha_0$
$E (u_\alpha)E(u_\alpha^*)$ is invertible in $\Mf$. Since $\Mf$ is finite, it follows 
that also $E (u_\alpha)E(u_\alpha^*)$ is invertible, which implies that $E (u_\alpha)$ 
is invertible. Then the unitary part 
$v_\alpha$ of $E(u_\alpha^*)$ ($E(u_\alpha^*)=
v_\alpha \left(E (u_\alpha)E(u_\alpha^*)\right)^{1/2}$) satisfies that $E (u_\alpha)v_\alpha \to 1$. 
Indeed, note that $v_\alpha =E (u_\alpha^*)
\left(E (u_\alpha) E(u_\alpha^*)\right)^{-1/2}$, and then $E (u_\alpha)v_\alpha=
E (u_\alpha)E(u_\alpha^*)(E (u_\alpha)E(u_\alpha^*))^{-1/2}\to 1$.
On the other hand, $[u_\alpha]\to [1]$ in ${\cal U}_\f=U_M/U_{\Mf}$ in the modular 
topology if and only if there exist unitaries 
$w_\alpha \in \Mf$ such that $\| u_\alpha w_\alpha -1\|_{E}\to 0$. 
Put $w_\alpha=v_\alpha$ as above,
then
$$
\begin{array}{rl}
\| u_\alpha v_\alpha -1\|_{E}^2 &=
\|E \left((v_\alpha^* u_\alpha^* -1)(u_\alpha v_\alpha -1)\right)\| \\
&
\le \|1-E(u_\alpha)v_\alpha\|+\|v_\alpha^* E(u_\alpha^*)-1\|
\end{array}
$$
which tend to zero, and therefore $\delta$ is a homeomorphism in
the modular topology of ${\cal U}_\f$.
\end{proof}
The following result will be useful in the study of these topologies
\begin{prop}
Let $u$ and $w$ be unitaries in $M$, then 
\begin{equation}\label{desigualdad}
\|\f \circ Ad(u^*) - \f \circ Ad(w^*)\| \le 2\|u - w\|_{E} \le 2\|u-w\|.
\end{equation}
\end{prop}
\begin{proof}{Proof.}
The second inequality is obvious. In order to prove the first note that
for any $x\in M$,
$$
|\f(u^*xu)-\f(w^*xw)|\le |\f(u^*x(u - w))|+ |\f((u^*-w^*)xw)|.
$$
Note that if $v$ is unitary, by the Cauchy-Schwarz inequality we have that  
$|\f(zv)| \le \f(zz^*)^{1/2}$ 
and $|\f(v^*z)|\le \f(z^*z)^{1/2}$. Applying these inequalities we obtain
$$
|\f(u^*x(u - w))|\le \f( (u^*-w^*)x^*x(u-v))^{1/2}=\f \circ E ( (u^*-w^*)x^*x(u-v))^{1/2},
$$
and
$$
|\f((u^*-w^*)xw)|\le \f \circ E ( (u^*-w^*)xx^*(u-w) )^{1/2}.
$$
Note that $(u^*-w^*)x^*x(u-v)\le \|x\|^2 (u^*-w^*)(u-v)$, and analogously for the other term. Thus we obtain
$$
|\f(u^*xu)-\f(w^*xw)|\le 2\|x\| \ \f\circ E ((u^*-w^*)(u-v))^{1/2}\le 2\|x\| \  \|E ((u^*-w^*)(u-v))\|^{1/2}.$$
\end{proof}
\begin{prop}
$\uf$ is complete in the usual topology (induced by the usual norm of $M$)
\end{prop}
\begin{proof}{Proof.}
One has continuous local cross sections 
$\sigma_{[u]} :{\cal V}_{[u]} \subset \uf \to U_M$ defined on a neighborhood 
${\cal V}_{[u]}$ of $[u] \in \uf$. 
If $[u_n]$ is a Cauchy net in $\uf$, choose $[u_{k_0}]$ such that for $n\ge k_0$,  
$[u_n] \in {\cal V}_{[u_{k_0}]}$. Then since $\sigma=\sigma_{[u_{k_0}]}$ is continuous, 
$\sigma ([u_n])$ is a Cauchy sequence in $U_M$ in the norm topology, therefore convergent 
to a unitary $v$ in $M$. Then clearly $[u_n]$ converges to $[v]$.
\end{proof}

Let us consider the same question for the modular topology. Of course, if the index of $E$ is 
finite, one obtains that $\uf$ is closed (and complete). In the general case, one expects to 
obtain other elements in the closure of $\uf$ with the modular topology.  Denote by $X_M$ the 
completion of the pre-Hilbert C$^*$-module $M$ to a Hilbert C$^*$-module. Note that $M$ acts 
on $X_M$ by left multiplication, and can be viewed as a closed subalgebra of the algebra of 
adjointable operators of $X_M$.

Suppose that $u_n$ is a sequence of unitaries converging to an element $x\in X_M$. Note that 
this implies that $<x,x>=1$, i.e. $x$ lies in the unit sphere of $X_M$. Put $\f_x \in M^*$, 
$\f_x (a)=\f(<x,ax>)$. Clearly $\f_x$ is a state of $M$.

\begin{prop} 
All elements in the closure of  $\uf$ with the modular topology are of the form $\f_x$ 
for some $x$ in the unit sphere of $X_M$. Such states $\f_x$ are normal. If $M$ is finite, 
then $\f_x$ is also faithful.
\end{prop}
\begin{proof}{Proof.}
An argument similar to the one in the previous proposition, shows that a Cauchy sequence $[u_n]$ 
for the modular topology yields another sequence
 of unitaries $u'_n$ in $U_M$ such that $[u_n]=[u'_n]$ and $u'_n$ form a Cauchy sequence in 
 $U_M$ for the norm $\| \ \|_{E}$ (this is clear using the continuous cross sections available 
 for this topology as well). Therefore
$u'_n$ converge to some element $x$ in the unit sphere of $X_M$. 
Now using the above result,
$$
\|\f \circ Ad(u_n^{' *}) -\f \circ Ad(u_k^{' *})\|\le 2\|u'_n -u'_k\|_{E},
$$
and therefore $\f\circ Ad(u_n^{' *})$ is a Cauchy sequence in $M_* \subset M^*$. Then 
$\f\circ Ad(u_n^{' *})$ converges to a normal state $\psi$.
Then $\f(u_n^{' *}au'_n)$ converges to  $\psi (a)$. On the other hand, 
$\f(u_n^{' *}au'_n)=\f \circ E(u_n^{' *}au'_n)=\f(<u'_n,au'_n>)$, which converges to 
$\f_x(a)$ by the continuity of the scalar product.

Suppose now that $M$ is finite, and fix a faithful tracial and normal state $\tau$. 
Pick $\f_x$ with $x$ a limit of unitaries $u_n$ as above. 
Then if $\f_x(a^*a)=0$, 
one has that $<x,xa^*a>=\linebreak
<ax,ax>=0$, which implies that $au_n$ tends to zero in the 
norm $\| \ \|_{E}$. In other words, $E(u_n^*a^*au_n) \to 0$. Note that $\tau \circ E=\tau$, and therefore
$\tau(E(u_n^*a^*au_n))=\tau (u_n^*a^*au_n)=\tau(a^*a)=0$, i.e. $a=0$.
\end{proof}

We do not know if these states in the closure of $\uf$ for the modular topology are in general 
faithful. There are other cases other than the finite case in which this happens. To prove our 
result we need the following lemma,
which was proven in \cite{acs}.

\begin{lemm}
The Jones projection $e$ associated to $E$ is finite in $M_1$.
\end{lemm}
If $E:M \to N \subset M$ is a conditional expectation, the normalizer ${\cal N}(E)$ is the group
$$
{\cal N}(E)=\{u \in U_M : Ad(u^*)\circ E \circ Ad(u)= E\}.
$$

\begin{prop}
If $\Mf$ and ${\cal N}(E)$ generate $M$ as a von Neumann algebra, then the states 
$\f_x$ in the closure of $\uf$ in the modular topology are faithful.
\end{prop}
\begin{proof}{Proof.}
Represent $M$ and $\Mf$ in $H_\f$ as in the basic construcion. As in the proposition above, 
one needs to show that if $u_n$ are unitaries in $M$ converging to $x \in X_M$ in the norm $\| \ \|_{E}$, 
and $a \in M$ such that $E(u_n^*a^*au_n)\to 0$, then it must be $a=0$.
We claim that under this hypothesis $au_n$ tends to zero in the strong operator topology. Note that 
 $au_ne$ tends to zero in norm, and therefore
$au_n be=au_ne b \to 0$ for all $b\in \Mf$. On the other hand,
if $v\in {\cal N}(E)$ then also $au_n ve \to 0$ in the norm of $M_1$. Indeed, 
$(au_n ve)^*au_n ve=e v^*u_n^*a^*au_nve=E(v^*u_n^*a^*au_nv)e=v^*E(u_n^*a^*au_n)ve \to 0$. 
Since $\Mf$ and ${\cal N}(E)$ generate $M$, it follows that $au_nxe$ tends to zero in 
norm for any $x\in M$. The claim is proven, using that the sequence is bounded in norm, 
and the fact that $Me$ is dense in $H_\f$. 
By the previous lemma, $e$ is finite, and therefore $(au_n)^*e \to 0$ in the strong operator 
topology (see \cite{kr}). Again using that the sequence is bounded, one has that 
$aa^*e=(au_n)(au_n)^*e \to 0$ strongly, that is $aa^*e=0$. Then $E(aa^*)e=e aa^*e=0$, 
which implies that $E(aa^*)=0$, and therefore $a=0$.
\end{proof}
There is an easy example of this situation. 
 \begin{exam}
Take $M=B(\ell^2(\zZ))$ and $E$ the conditional expectation onto the subalgebra of 
diagonal matrices in the canonical basis of $\ell^2(\zZ)$. This subalgebra is the 
centralizer of a state $\f$, 
which can be constructed by means of a density operator $a$ chosen as a diagonal trace 
class positive matrix, with different non zero entries in the diagonal, and with trace $1$. 
The bilateral shift and its integer powers normalize this expectation, and it is 
straightforward to verify that the diagonal matrices together with the powers of the 
bilateral shift generate $B(\ell^2(\zZ))$.
\end{exam}	
The inequality \ref{desigualdad} implies that one can embed $\uf$, 
both with the usual and the modular topologies, 
in the state space  of $M$ (with the norm topology of $M^*$).
However the usual and the modular topology of $\uf$ do not coincide 
with the norm topology of the state space. This fact is clear  in the following example. 
\begin{exam}
Let $M$ and $\f=Tr(a \ .)$ as in the preceeding example.
Let $t=I \oplus \sigma$ act on $\ell^2(\zZ)=\ell^2(\zN) \oplus \ell^2(\zN)$, where $\sigma$ 
is the unilateral shift. Denote by $q_n$ the $n \times n$ Jordan nilpotent, and $u_n$ the 
unitary operator on $\ell^2(\zN)$ having the unitary matrix $q_n+q_n^{* \ n-1}$ on the first 
$n\times n$ corner, and the identity matrix afterwards. Finally put $w_n =1\oplus u_n \in B(\ell^2(\zZ))$.
It is straightforward to verify that $w_n a w_n^*\to tat^*$ in the trace norm of 
$B(\ell^2(\zZ))$. This means that $\f \circ Ad(w_n^*)\to \f_t$ in the norm of $B(\ell^2(\zZ))^*$. 
But by the proposition above, it is clear that $\f_t$ does not belong to the closure of $\uf$ either 
in the modular or usual topology, since it is not faithful ($\ker t$ is not trivial).
\end{exam}

\vskip0.5cm

{\sit 
\noindent
Esteban Andruchow and Alejandro Varela\\
Instituto de Ciencias \\
Universidad Nacional de Gral. Sarmiento \\
J. A. Roca 850 \\
(1663) San Miguel \\
Argentina  \\
e-mail:  eandruch@ungs.edu.ar, avarela@ungs.edu.ar
}

\end{document}